\documentclass[a4paper,11pt]{article}
\usepackage{graphicx,amssymb,amsmath,amsthm,enumerate}
\usepackage[numbers,sort&compress]{natbib}
\usepackage[lmargin=31mm,rmargin=31mm,tmargin=34.5mm,bmargin=34.5mm]{geometry}
\renewcommand{\baselinestretch}{1.169}
\newtheorem{theorem}{Theorem}
\newtheorem{lemma}[theorem]{Lemma}
\newtheorem{cor}[theorem]{Corollary}
\newtheorem{prop}[theorem]{Proposition}

\newtheorem{conj}{Conjecture}
\newtheorem{open}{Open Problem}

\newcommand{\PP}[1]{\ensuremath{A\substack{\longrightarrow\\#1}B}}

\begin{document}
\author{
Ruy Fabila-Monroy\footnotemark[1] \and
David Flores-Pe\~naloza \footnotemark[1] \and
Clemens Huemer \footnotemark[3] \and
Ferran Hurtado \footnotemark[2] \and
Jorge Urrutia  \footnotemark[1] \and
David R.~Wood \footnotemark[4]}

\renewcommand{\thefootnote}{\fnsymbol{footnote}}

\footnotetext[1]{Instituto de Matem\'aticas, Universidad Nacional Aut\'onoma de M\'exico (\texttt{ruy@ciencias.unam.mx}, \texttt{\{dflores,urrutia\}@math.unam.mx}). }

\footnotetext[2]{Departament de Matem\`atica Aplicada II, Universitat
  Polit\`ecnica de Catalunya, Barcelona, Spain
  (\texttt{ferran.hurtado@upc.edu}). Partially supported by projects
  MEC MTM2009-07242 and Gen.\  Cat.\  DGR 2009SGR1040.}

\footnotetext[3]{Departament de Matem\`atica Aplicada IV, Universitat
  Polit\`ecnica de Catalunya, Barcelona, Spain
  (\texttt{clemens.huemer@upc.edu}). Partially supported by projects
  MEC MTM2009-07242 and Gen.\  Cat.\ DGR 2009SGR1040.}

\footnotetext[4]{Department of Mathematics and Statistics, The
  University of Melbourne, Melbourne, Australia
  (\texttt{woodd@unimelb.edu.au}). Supported by a QEII Research
  Fellowship from the Australian Research Council.}

\renewcommand{\thefootnote}{\arabic{footnote}}

\title{\bf Token Graphs}
\maketitle

\begin{abstract}
For a graph $G$ and integer $k\geq1$, we define the token graph $F_k(G)$ to be the graph with vertex set all $k$-subsets of $V(G)$, where two vertices are adjacent in $F_k(G)$ whenever their symmetric difference is a pair of adjacent vertices in $G$. Thus vertices of $F_k(G)$ correspond to configurations of $k$ indistinguishable tokens placed at distinct vertices of $G$, where two configurations are adjacent whenever one configuration can be reached from the other by moving one token along an edge from its current position to an unoccupied vertex. This paper introduces token graphs and studies some of their properties including: connectivity, diameter, cliques, chromatic number, Hamiltonian paths, and Cartesian products of token graphs.
\end{abstract}

\section{Introduction} \label{sec:intro}

Many problems in mathematics and computer science are modeled by moving objects on the vertices of a graph according to certain prescribed rules. In ``graph pebbling'', a pebbling step consists of removing two pebbles from a vertex and placing one pebble on an adjacent vertex; see \cite{surpeb} and \citep{surpebnew} for surveys. Related pebbling games have been used to study rigidity \citep{pebblealg,pebblejacobs}, motion planning \citep{auletta99linear,papa94}, and as models of computation \citep{models}. 
In the ``chip firing game'', a vertex $v$ fires by distributing one chip to each of its neighbors (assuming the number of chips at $v$ is at least its degree).
This model has connections with matroids, the Tutte polynomial, and  mathematical physics; see \citep{merino} for a survey. 

\begin{figure}[!h] 
\label{fig:TokenPath}
  \begin{center}
   \includegraphics{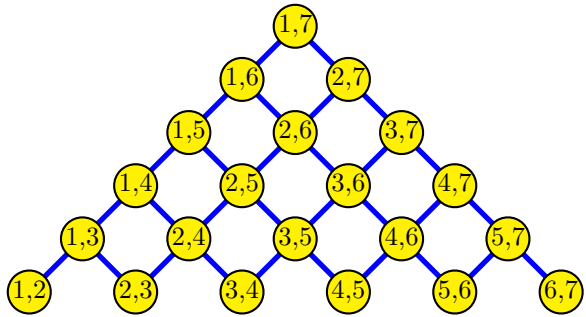}
\vspace*{-1ex}
   \caption{The $2$-token graph of the 7-vertex path.}
  \end{center}
\end{figure}

In this paper we study a model in which $k$ indistinguishable tokens
move from vertex to vertex along the edges of a graph. This idea is
formalized as follows. For a graph\footnote{We consider undirected,
  simple and finite graphs; see \citep{diestel}. A \emph{$k$-set} is a
  set with cardinality $k$. For a set $S$, let $\binom{S}{k}$ be the
  set of all $k$-sets contained in $S$. Let $[a,b]:=\{a,a+1,\dots,b\}$
  and $[n]:=[1,n]$. For sets $A$ and $B$, let $A \triangle B:=(A \cup
  B)\setminus(A \cap B)$. }  $G$ and integer $k\geq1$, we define
$F_k(G)$ to be the graph with vertex set $\binom{V(G)}{k}$, where two
vertices $A$ and $B$ of $F_k(G)$ are adjacent whenever their symmetric difference $A \triangle B$ is a pair $\{a,b\}$ such that $a \in A$, $b \in B$  and $ab\in E(G)$.
Thus the vertices of $F_k(G)$ correspond to configurations of $k$ indistinguishable tokens placed at distinct vertices of $G$, where two configurations are adjacent whenever one configuration can be reached from the other by moving one token along an edge from its current position to an unoccupied vertex. We thus call $F_k(G)$ the {\em $k$-token graph of $G$}. See Figure~\ref{fig:TokenPath} for an example.

The aim of this paper is to introduce token graphs and study some of their properties. 
We make the following contributions:
We prove tight lower and upper bounds on the diameter of token graphs (Section~\ref{sect:con}).
We prove tight lower bounds on the connectivity of token graphs (Section~\ref{sect:con}). 
We characterize the cliques in token graphs in terms of the cliques in the original graph, and derive an exact formula for the clique-number of a token graph (Section~\ref{sect:clique}).
We present upper and lower bounds on the chromatic number of token graphs, and conclude that every token graph has chromatic number at least (roughly) half the chromatic number of the original graph and
at most the chromatic number of the original graph (Section~\ref{sect:chrom}). 
We establish sufficient conditions for the existence or non-existence of a Hamiltonian path in various token graphs (Section~\ref{sect:ham}).
We show that token graphs contain certain Cartesian products as induced subgraphs (Section~\ref{sect:prod}).
Finally, we suggest some new research problems, mostly related to graph reconstruction (Section~\ref{sect:end}).

A key example in our study is when $G$ is a complete graph. 
Then the token graph is called a Johnson graph, which is widely studied due to connections with coding theory.
The Johnson graph $J(n,k)$ is the graph whose vertices are the $k$-subsets of an $n$-set, where
two vertices $A$ and $B$ are adjacent whenever $|A\cap B|=k-1$ (or alternatively, if $|A\triangle B|=2$). 
Observe that $F_k(K_n) \simeq J(n,k)$. Many results in this paper generalize known properties of Johnson graphs.

\section{Basic Properties} \label{sect:prop} 

Throughout this paper, $G$ is a graph with $n$ vertices and $k$ is a positive integer. 
To avoid trivial cases, we assume that $n\geq k+1$. 
The number of vertices in $F_k(G)$ is:
$$|V(F_k(G))|=\binom{n}{k}\enspace.$$

To calculate the number of edges in $F_k(G)$, 
charge each edge $AB$ of $F_k(G)$ to the unique edge $ab$ of $G$, for which $A\triangle B=\{a,b\}$.
The number of edges of $F_k(G)$ charged to $ab$ is $\binom{n-2}{k-1}$.
Hence
$$|E(F_k(G))|=\binom{n-2}{k-1}|E(G)|\enspace.$$

The neighborhood of each vertex $A$ of $F_k(G)$ is 
$$\{A\setminus\{v\}\cup\{w\}:v\in A,w\in V(G)\setminus A,vw\in E(G)\}\enspace.$$
Thus the degree of $A$ in $F_k(G)$ equals the number of edges between $A$ and $V(G)\setminus A$. 
Straightforward bounds on the minimum and maximum degree of $F_k(G)$ follow.

With only one token, the resulting token graph is isomorphic to $G$. Thus
\begin{equation}
\label{eqn:1token}
F_1(G) \simeq G\enspace.
\end{equation}

Since two vertices $A$ and $B$ are adjacent in $F_k(G)$ if and only if $V(G)\setminus A$ and $V(G)\setminus B$ are adjacent 
in $F_{n-k}(G)$, 
\begin{equation} \label{eqn:k,n-k}
F_k(G) \simeq F_{n-k}(G)\enspace.
\end{equation}




We sometimes use \eqref{eqn:k,n-k} to assume that $k \le\frac{n}{2}$. Also note that \eqref{eqn:1token} and \eqref{eqn:k,n-k} imply two known properties of the Johnson graph: $J(n,1)\simeq K_n$ and $J(n,k) \simeq J(n,n-k)$.

At times, we study the token graph that arises when tokens are fixed at certain vertices. 
Given a set $X \subseteq V(G)$ with $|X|=r\le k$, we define $F_k(G,X)$ to be the
subgraph of $F_k(G)$ induced by the vertices of $F_k(G)$ that contain $X$ as a subset.  
This definition can be interpreted as having $r$ tokens fixed at $X$, and $k-r$ tokens moving on $G-X$.
Hence
\begin{equation} \label{eqn:fixed}
F_k(G,X) \simeq F_{k-r}(G-X)\enspace.
\end{equation}

\section{Connectivity and Diameter}
\label{sect:con}

In this section we establish tight bounds on the connectivity and diameter of $F_k(G)$ in terms of the same parameters in $G$. 

The following notation will be helpful. 
Let $A$ be a $k$-set in a graph $G$. 
Let $P$ be an $ab$-path in $G$ such that $a\in A$ and $b\not\in A$.
Let $A':=A\setminus\{a\}\cup\{b\}$.
Say $A\cap P=\{v_1,\dots,v_q\}$ ordered by $P$ (although not
necessarily consecutive in $P$), where $v_1=a$. 
Let $A\substack{\longrightarrow\\P}A'$ be the path between $A$ and $A'$ in $F_k(G)$ corresponding to the following sequence of token moves:
First move the token at $v_q$ along $P$ to $b$, 
then for $i=q-1,q-2,\dots,1$ 
move the token at $v_i$ along $P$ to $v_{i+1}$.
Each move is along a path containing no tokens. 
Thus these moves correspond to a path in $F_k(G)$. 
Observe that this path terminates at $A'$. 
Each edge in  $A\substack{\longrightarrow\\P}A'$ corresponds to an edge in $P$.
Thus the length of $A\substack{\longrightarrow\\P}A'$ equals the length of $P$. 

\begin{theorem}\label{thm:Diameter}
Let $G$ be a connected graph with diameter $\delta$. 
Then $F_k(G)$ is connected with diameter at least $k(\delta-k+1)$ and at most $k\delta$. 
\end{theorem}

\begin{proof}
We prove the upper bound by induction on $|A\triangle B|$ with the following hypothesis: 
``for all vertices $A,B$ of $F_k(G)$ there is an $AB$-path in $F_k(G)$ of length at most $\frac{\delta}{2}|A\triangle B|$.''
This implies that $F_k(G)$ is connected with diameter at most $k\delta$. 

If $A\triangle B=\varnothing$ then $A=B$ and there is nothing to prove.
Now assume that $A\triangle B\neq\varnothing$.
Since $G$ is connected there is a path $P$ between some vertex $a\in A-B$ and $b\in B-A$. 
Thus $A\substack{\longrightarrow\\P}A'$ is a path in $F_k(G)$ from $A$ to  $A':=A\setminus\{a\}\cup\{b\}$.
Observe that $|A'\triangle B|=|A\triangle B|-2$.
By induction there is path between $A'$ and $B$ in $F_k(G)$ of length at most $\frac{\delta}{2}|A'\triangle B|=\frac{\delta}{2}|A\triangle B|-\delta$. 
Since the length of $A\substack{\longrightarrow\\P}A'$ equals the length of $P$, which is at most $\delta$, 
there is path between $A$ and $B$ in $F_k(G)$ of length at most $\frac{\delta}{2}|A\triangle B|$.

Now we prove the lower bound. 
Let $x$ and $y$ be vertices at distance $\delta$ in $G$.
For $i\in[0,\delta]$, let $V_i$ be the set of vertices in $G$ at distance $i$ from $x$. 
Thus $V_0=\{x \} $ and $y \in V_\delta$. 
Let $d(v)$ be the distance between $x$ and each vertex $v$. 

Let $a$ be the minimum index such that $|V_0\cup\dots\cup V_a| \ge k$.
Likewise, let $b$ be the maximum index such that $|V_b\cup\dots\cup V_\delta| \ge k$. 
Let $A$ be a subset of $V_0\cup\dots\cup V_a$ such that $V_0\cup\dots\cup V_{a-1}\subset A$.
Let $B$ be a subset of $V_b\cup\dots\cup V_\delta$ such that $V_{b+1}\cup\dots\cup V_{\delta}\subset B$.



Consider any path from $A$ to $B$ in $F_k(G)$. 
Each token initially at a vertex $v\in A$ is moved to some vertex $v'\in B$.
Since edges in $G$ are either within some set $V_i$ or between sets $V_i$ and $V_{i+1}$, 
at least $d(v')-d(v)$ moves are needed to move the token from $v$ to $v'$. 
Thus $F_k(G)$ has diameter at least 
$$\sum_{v \in A} (d(v')-d(v))=\sum_{w \in B}d(w)-\sum_{v \in A} d(v)\enspace.$$
The first summation is minimized when $b=\delta-k+1$ and $|V_j|=1$ for all $j \geq b$.
The second summation is maximized when $a=k-1$ and $|V_i|=1$ for all $i\leq a$. 
Thus $F_k(G)$ has diameter at least 
$$\sum_{j=\delta-k+1}^{\delta} \!\!\!\!j\;-\;\sum_{i=0}^{k-1} i=k(\delta-k+1)\enspace.$$
\end{proof}

\begin{figure}[!h] 
\label{diamexample}
  \begin{center}
   \includegraphics{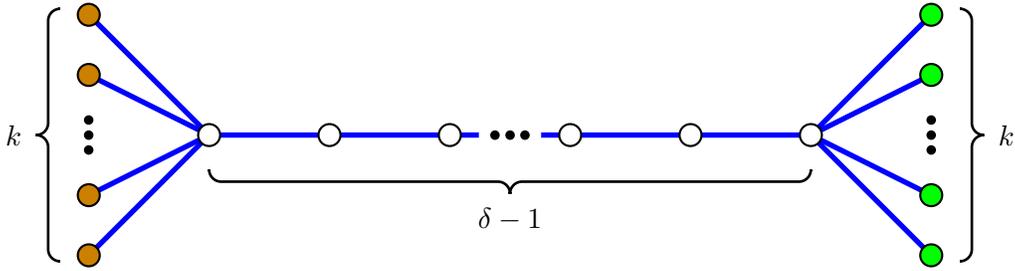}
   \caption{Configurations at distance $\delta k$ in $F_k(T)$.}
  \end{center}
\end{figure}

Note that both bounds in Theorem~\ref{thm:Diameter} are achievable. If
$P_{\delta+1}$ is the path on $\delta+1$ vertices and $k\leq
\delta+1$, then $P_{\delta+1}$ has diameter $\delta$ and
$F_k(P_{\delta+1})$ has diameter $k(\delta-k+1)$. And, as illustrated
in Figure~\ref{diamexample}, if $T$ is the tree obtained by adding $k$
vertices adjacent to each endpoint of $P_{\delta-1}$, then $T$ has
diameter $\delta$ and $F_k(T)$ has diameter $\delta k$. 


We now consider the connectivity of $F_k(G)$ when $G$ is highly connected. 

\begin{lemma}
\label{lem:TwoDisjointPaths}
Let $A$ be a $k$-set in a graph $G$. 
Let $a$ and $b$ be vertices of $G$  such that $a\in A$ and $b\not\in A$.
Let $P$ and $Q$ be internally disjoint $ab$-paths in $G$.
Then $A\substack{\longrightarrow\\P}A'$
and $A\substack{\longrightarrow\\Q}A'$ are internally disjoint paths in $F_k(G)$, where $A'=A\setminus\{a\}\cup\{b\}$. 
\end{lemma}

\begin{proof}
First suppose that $|P\cap A|\geq2$. 
Let $P\cap A=\{v_1,\dots,v_p\}$ ordered by $P$, where $v_1=a$.
Consider an arbitrary internal vertex $R$ of $A\substack{\longrightarrow\\P}A'$.
Then $|R\cap P|=p$ and $R$ contains a vertex in the sub-path $(v_p,b]$ of $P$.
Thus $R$ does not contain $\{v_2,\dots,v_p\}$.
On the other hand, $\{v_2,\dots,v_p\}$ is contained in every vertex of $A\substack{\longrightarrow\\Q}A'$.
Thus $A\substack{\longrightarrow\\P}A'$ and $A\substack{\longrightarrow\\Q}A'$ are internally disjoint, and we are done.
Similarly if $|Q\cap A|\geq2$. 

Now assume that $P\cap A =Q\cap A=\{a\}$.
Without loss of generality, $P$ is not the edge $ab$.
Thus $P\setminus\{a,b\}\neq\varnothing$. 
Thus every internal vertex of $A\substack{\longrightarrow\\P}A'$ contains some vertex in $P\setminus\{a,b\}$.
On the other hand, no internal vertex of $A\substack{\longrightarrow\\Q}A'$ contains a vertex in $P\setminus\{a,b\}$.
Thus $A\substack{\longrightarrow\\P}A'$ and $A\substack{\longrightarrow\\Q}A'$ are internally disjoint.
\end{proof}

We need the following technical result in the proof of
Lemma~\ref{lem:k-1} below.

\begin{lemma}
\label{lem:RedBlue}
Let $H$ be a complete bipartite graph with colour classes $Y$
and $Z$, where $|Y|<|Z|$. 
Suppose that the edges of $H$ are coloured red and blue, such
that each vertex in $Y$ is incident to at most one red edge. 
Then $H$ contains a set $M$ of blue edges, such that
each vertex in $Y$ is incident to exactly one edge in $M$, and
the union of the red edges and $M$ is acyclic.
\end{lemma}

\begin{proof}
We proceed by induction on $|Y|$. 
The base case is trivial. Since there are more vertices in $Z$ than
red edges, some vertex $x\in Z$ is incident to no red edge.
Let $v$ be any vertex in $Y$. 
Let $vw$ be the red edge incident to $v$ (if any). 
Let $H':=(H-v)-x$. 
Let $R$ and $R'$ be the sets of red edges in $H$ and $H'$ respectively.
By induction, there is a set $M'$ of blue edges in $H'$, such
that each vertex in $Y-v$ is incident to exactly one edge in $M'$, and
$R'\cup M'$ is acyclic.
Let $M:=M'\cup\{vx\}$. Thus $v$ (and every vertex in $Y$) is incident
to exactly one edge in $M$. 
Since $x$ is incident to no red edge, $M\cup R$ is obtained from
$M'\cup R'$ by adding the edges $xv$ and $vw$ (if it exists). 
Thus $M\cup R$ is acyclic.
\end{proof}

A \emph{chord} of a path $P$ in a graph $G$ is an edge
$vw\in E(G)\setminus E(P)$, such that
both $v$ and $w$ are in $P$, but the endpoints of $P$ are not $v$ and
$w$. Thus $P$ is \emph{chordless} if the subgraph of $G$ induced
by $V(P)$ has maximum degree at most $2$. 

\begin{lemma}
\label{lem:k-1} 
Let $G$ be a $t$-connected graph. 
Let $A$ and $B$ be vertices of $F_k(G)$ such that $|A\triangle B|=2$. 
Then there are $t$ internally disjoint $AB$-paths in $F_k(G)$.
Moreover, if $t\geq k$ then there are $k(t-k+1)$ internally disjoint $AB$-paths in $F_k(G)$.
\end{lemma}

\begin{proof}
Let $a$ and $b$ be the vertices in $A\setminus B$ and $B\setminus A$
respectively. 
By Menger's Theorem, $G$  contains internally disjoint $ab$-paths $P_1,\dots,P_t$. 
Thus $\PP{P_1},\dots,A\substack{\longrightarrow\\P_t}B$ are internally disjoint $AB$-paths in $F_k(G)$ by Lemma~\ref{lem:TwoDisjointPaths}. This proves the first claim. 

Now assume that $t\geq k+1$. 
As illustrated in Figure~\ref{fig:Proof}, let
$P_1,\dots,P_s,Q_1,\dots,Q_\ell$ be a set of internally disjoint $ab$-paths, such that $s+\ell\geq t$, where each of the paths $P_1,\dots,P_s$ do not
intersect $A\cap B$, and each of the paths $Q_1,\dots,Q_\ell$ 
do intersect $A\cap B$. 
Choose such a set of paths such that $s+\ell$ is maximal and each path is chordless. 

\begin{figure}[!h] 
 \begin{center}
   \includegraphics{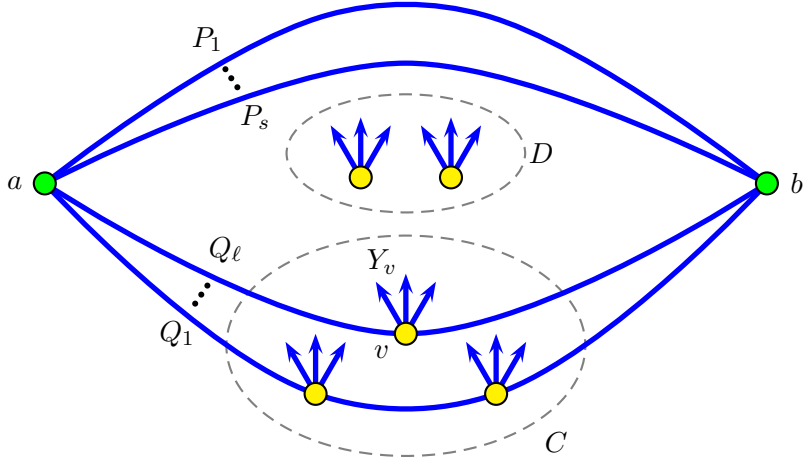}
  \caption{Illustration for the proof of Lemma~\ref{lem:k-1}.}
\label{fig:Proof}
  \end{center}
\end{figure}

Let $C$ be the set of vertices in $A\cap B$
that intersect one of  $Q_1\cup\dots\cup Q_\ell$. 
Since $Q_1,\dots,Q_\ell$ are internally disjoint, each vertex in $C$
is in exactly one of $Q_1,\dots,Q_\ell$.
Let $D$ be the set of vertices in $A\cap B$
that do not intersect $Q_1\cup\dots\cup Q_\ell$. Thus $C$ and $D$ partition $A\cap
B$. Hence $|C|+|D|=k-1$ and $\ell\leq|C|\leq
k-1$ and $$s\geq t-\ell\geq t-|C|=t-(k-1-|D|)=t-k+1+|D|\enspace.$$


The $AB$-paths that we construct
in $F_k(G)$ are of three types. The first and second types are
straightforward. By
Lemma~\ref{lem:TwoDisjointPaths}, 
$$A\substack{\longrightarrow\\P_1}B,\dots,A\substack{\longrightarrow\\P_s}B,
A\substack{\longrightarrow\\Q_1}B,\dots,A\substack{\longrightarrow\\Q_\ell}B$$
are internally disjoint $AB$-paths in $F_k(G)$, called \emph{type-$P$}
and \emph{type-$Q$} paths respectively.
Note that since $P_i$ avoids $A\cap B$, the path
$A\substack{\longrightarrow\\P_i}B$ in $F_k(G)$ corresponds to the sequence of token configurations obtained by simply moving the token from $a$ along $P_i$ to $b$. 
 
For each vertex $v\in A\cap B$, we construct a set of \emph{type-$R$} paths in
$F_k(G)$ between $A$ and $B$ as follows. 

First consider $v\in C$. Then $v\in Q_i$ for exactly one
$i\in[\ell]$. 
Let $Y_v:=N_G(v)\setminus(A\cap B)\setminus Q_i$. 
Since $G$ is $t$-connected, $|N_G(v)|=\deg_G(v)\geq t$.
Since $Q_i$ is chordless, $v$ has only two
neighbours in $Q_i$. 
Since $|A\cap B|=k-1$ and $v\in (A\cap B)-N_G(v)$, we have 
$|Y_v|\geq t-(k-2)-2=t-k$. 

Now consider $v\in D$.
Let $Y_v:=N_G(v)\setminus(A\cup B)$. 
Since $G$ is $t$-connected, $|N_G(v)|=\deg_G(v)\geq t$.
We have $|A\cup B|=k+1$ and $v\in (A\cup B)\setminus N_G(v)$.
Moreover, since  $s+\ell$ is maximal, the path $(a,v,b)$ is not in
$G$. 
Thus $a\not\in N_G(v)$ or $b\not\in N_G(v)$. 
Hence $|Y_v|\geq t-(k-1)=t-k+1$. 

Now let $Y_v$ be a subset of itself with exactly 
$t-k$ vertices if $v\in C$, and exactly $t-k+1$ vertices if $v\in D$. 
Note that $Y_v\neq\varnothing$ since $t\geq k+1$. 
And by construction, $a\not\in Y_v$ and
$b\not\in Y_v$. 

Let $H_v$ be the complete bipartite graph with colour classes $Y_v$
and $[s]$. Colour the edges of $H_v$ as follow. 
If $y\in P_i$ for some $y\in Y_v$ and $i\in[s]$, then  colour the edge
$yi$ in $H_v$ red. 
Colour every other edge in $H_v$ blue. 
Since each vertex in $Y_v$ is in at most one of $P_1,\dots,P_s$, each
vertex in $Y_v$ is incident to at most one red edge in $H_v$. 
We now verify that Lemma~\ref{lem:RedBlue} is applicable to $H_v$ with
$Y=Y_v$ and $Z=[s]$. 
If $v\in C$ then $|Z|=s\geq t-k+1>t-k=|Y_v|$. 
If $v\in D$ then $|D|\geq1$ and 
$|Z|=s\geq t-k+1+|D| >t-k+1=|Y_v|$. 
In both cases, Lemma~\ref{lem:RedBlue} is applicable . 

Thus there is set $M_v$ of blue edges in $H_v$, such that 
each vertex in $Y_v$ is incident to exactly one edge in $M_v$, and
the union of the red edges and $M$ is an acyclic subgraph of $H_v$. 
Note that $|M_v|=|Y_v|$. 
For each edge $yi\in M_v$, let $R\langle{v,y\rangle}$ be
the \emph{type-$R$} path in $F_k(G)$ corresponding to the following
token moves (where all the tokens at $(A\cap B)\setminus\{v\}$ are stationary):

(1) move the token at $v$ to $y$,

(2) move the token at $a$ along the path $P_i$ to $b$,

(3) move the token at $y$ back to $v$. 

We now prove that the type-$R$ paths are internally disjoint. 
Suppose to the contrary that $R\langle{v,y\rangle}$ and
$R\langle{v',y'\rangle}$ 
share a common internal vertex, for some $(v,y)\neq(v',y')$.
Thus $yi$ is an edge of $M_v$, and $y'i'$ is an edge in $M_{v'}$, for
some $i,i'\in[s]$.
Each internal vertex in $R\langle{v,y\rangle}$ consists of $(A\cap
B)\setminus\{v\}\cup\{y\}$ plus some vertex in $P_i$.
Hence $$(A\cap B)\setminus\{v\}\cup\{y,x\}=
(A\cap B)\setminus\{v'\}\cup\{y',x'\}$$
for some vertices $x$ in $P_i$ and $x'\in P_{i'}$.  
Since $A\cap B\cap Y_{v'}=\varnothing$ and $y'\in Y_{v'}$, we have
$y'\not\in (A\cap B)\setminus\{v\}$, implying $y'\in\{x,y\}$.
Since $A\cap B\cap P_{i'}=\varnothing$ and $x'\in P_{i'}$, we have
$x'\not\in (A\cap B)\setminus\{v\}$, implying $x'\in\{x,y\}$.
Thus $\{x',y'\}=\{x,y\}$, implying $(A\cap B)\setminus\{v\}=(A\cap
B)\setminus\{v'\}$ and $v=v'$. 
Hence $yi$ and $y'i'$ are edges in the same set $M_v$. 
Since each vertex in $Y_v$ is incident to exactly one edge in $M_v$, 
we have $y\neq y'$. 
Thus $x=y'$ and $y=x'$, implying $y\in P_{i'}$ and $y'\in P_i$.
Hence, in $H_v$, the edges $yi'$ and $y'i$ are both red.
Since $yi$ and $y'i'$ are blue edges, $i\neq i'$.
Thus $(y,i,y',i')$ is a blue--red--blue--red cycle in $H_v$ with both blue edges in $M_v$. 
This contradiction proves that the type-$R$ paths are pairwise disjoint. 


We now prove that each type-$P$ path is internally disjoint from
each type-$R$ path. 
Suppose on the contrary that some path \PP{P_i} intersects some path
$R\langle{v,y}\rangle$ at an internal vertex in common.
Now $v$ is in every internal vertex of \PP{P_i}
(that is, the token at $v$ never moves in this sequence). 
On the other hand, $v$ is in no internal vertex of
$R\langle{v,y}\rangle$.  This contradiction proves that each type-$P$ path is internally disjoint from
each type-$R$ path. 

We now prove that each type-$Q$ path is internally disjoint from
each type-$R$ path. Suppose on the contrary that some path
$A\substack{\longrightarrow\\Q_i}B$ 
intersects some path  $R\langle{v,y}\rangle$ at an internal vertex $W$ in common.
Let $yj$ be the blue edge in $M_v$, where $j\in[s]$. 
By the construction of $H_v$,  we have $y\not\in P_j$. 
If $v\not\in Q_i$ then $v$ is in every internal vertex of \PP{Q_i}
(that is, the token at $v$ never moves in this sequence). 
On the other hand, $v$ is in no internal vertex of
$R\langle{v,y}\rangle$. 
Thus $v\in Q_i$ and $v\in C$. 
Since $y\in Y_v$ and $Y_v\cap Q_i=\varnothing$, we have $y\not\in
Q_i$. 
Similarly, since $y\in Y_v$ and $Y_v\cap A\cap B=\varnothing$, we have $y\not\in
A\cap B$. 
Every internal vertex of $R\langle{v,y}\rangle$ contains $y$.
But every internal vertex of \PP{Q_i} is contained in $(A\cap B)\cup
Q_i$. This contradiction proves that each type-$Q$ path is internally disjoint from each
type-$R$ path.

We have $s$ type-$P$ paths and $\ell$ type-$Q$ paths. 
For each $v\in C$ we have $t-k$ type-$R$ paths, and 
for each $v\in D$ we have $t-k+1$ type-$R$ paths.
In total, the number of $AB$-paths in $F_k(G)$ is 
\begin{align*}
&\;s+\ell+|C|(t-k)+|D|(t-k+1)\\
=&\;s+\ell+(|C|+|D|)(t-k)+|D|\\
=&\;s+\ell+(k-1)(t-k)+|D|\\
\geq &\; t+(k-1)(t-k)\\
=&\; k(t-k+1)\enspace.
\end{align*}
Therefore we have $k(t-k+1)$ pairwise internally disjoint $AB$-paths in $P_k(G)$.
\end{proof}

\begin{theorem} 
\label{thm:LowConnectivity}
If $G$ is $t$-connected, then $F_k (G)$ is $t$-connected for all $k\geq1$. 
\end{theorem}

\begin{proof}
By \eqref{eqn:k,n-k}, we may assume that  $k\leq\frac{n}{2}$.
Let $\mathcal{C}$ be a minimum (vertex) cut-set of $F_k(G)$. 
We will prove that $|\mathcal{C}| \ge t$, implying $F_k (G)$ is $t$-connected.

Let $A$ and $B$ be vertices in distinct components of $F_k(G)-\mathcal{C}$, such that $|A\triangle B|$ is minimum. 
If $|A \triangle B|=2$ then by Lemma~\ref{lem:k-1}, 
there are $t$ internally disjoint $AB$-paths in $F_k(G)$, implying $|\mathcal{C}|\geq t$. 
Now assume that $|A \triangle B|=2r\geq 4$.

Let $A \setminus B=\{a_1, \dots ,a_r\}$ and $B \setminus A=\{b_1, \dots ,b_r\}$.
For each $i\in[r]$ and each vertex $x\in V(G)\setminus(A \cup B)$, 
define   $A_{i,x}:=A \setminus \{a_i\} \cup \{x\}$ and     $B_{i,x}:=B \setminus \{b_i\} \cup \{x\}$. 
Suppose that  $A_{i,x}\notin \mathcal{C}$ and $B_{j,x} \notin \mathcal{C}$ for some $i,j$ and $x$.
Since $A \triangle A_{i,x}=\{a_i,x\}$, by the minimality of $|A\triangle B|$, 
$A$ and $A_{i,x}$ are in  the same component  of $F_k(G) - \mathcal{C}$. 
Similarly, $B$ and $B_{j,x}$ are in the same component of $F_k(G)-\mathcal{C}$. 
Also $A_{i,x} \triangle B_{j,x}=(A\triangle B) \setminus\{a_i,b_j\}$, implying $|A_{i,x} \triangle B_{j,x}|=2(r-1)$.
Thus $A_{i,x}$ and $B_{j,x}$ are in the same component of $F_k(G)$.
Hence $A$ and $B$ are in the same component of $F_k(G)$. 
This contradiction implies that at least one of $A_{i,x}$ and
$B_{j,x}$ is in $\mathcal{C}$ for all $i,j$ and $x$.

Hence, for each $x\in V(G)\setminus(A\cup B)$, the set $\mathcal{C}$
includes all of $\{A_{i,x}:i\in[r]\}$ or all of
$\{B_{j,x}:j\in[r]\}$. 
This gives $r(n-k-r)$ vertices in $\mathcal{C}$ (since $|A\cup B|=k+r$). 

Now for each $i\in[r]$ and $j\in[r]$, let $A_{i,j}:=A \setminus \{a_i\} \cup \{b_j\}$ and     $B_{i,j}:=B \setminus \{b_i\} \cup \{a_j\}$.  
There are $2r^2$ such sets.     
Suppose that $A_{i,j}\notin \mathcal{C}$.
Thus $A$ and $A_{i,j}$ are  in the same component of $F_k(G)-\mathcal{C}$    since $A\triangle A_{i,j}=\{a_i,b_j\}$.     
And $B$ and $A_{i,j}$ are in the same component of     $F_k(G)-\mathcal{C}$ since $|B \triangle A_{i,j}|=2(r-1)$.
Hence $A$ and $B$ are in the  same component of     $F_k(G)-\mathcal{C}$. 
This contradiction proves that $A_{i,j}\in \mathcal{C}$.
Similarly, $B_{i,j}\in \mathcal{C}$. 
These vertices were not counted in $\mathcal{C}$ earlier. 
Thus, since $r\geq2$ and $k\leq\frac{n}{2}$,
$$|\mathcal{C}|
\geq r(n-k-r)+2r^2
> r(n-k)
\geq 2(n-k)
\geq n> t\enspace.$$ 
Therefore $F_k(G)$ is $t$-connected.
\end{proof}

Theorem~\ref{thm:LowConnectivity} is best possible when $t\leq k$.
Let $G$ be a $t$-connected graph containing an edge cut $S$ of $t$ edges, 
such that the union $A$ of some components of $G-S$ has exactly $k$ vertices
(for example, take a matching of $t$ edges between two disjoint copies of $K_k$). 
Then $A$ has degree $t$ in $F_k(G)$. Thus $F_k(G)$ has connectivity exactly $t$.
We now prove a stronger bound for large $t$ and sufficiently large graphs.

\begin{theorem} 
\label{thm:BigConnectivity}
If $G$ is $t$-connected and $t\geq k$ and $n\geq \tfrac{1}{2}kt$, then $F_k (G)$ is $k(t-k+1)$-connected.
\end{theorem}

\begin{proof}
Let $\mathcal{C}$ be a minimum (vertex) cut-set of $F_k(G)$. 
Let $A$ and $B$ be vertices in distinct components of $F_k(G)-\mathcal{C}$, such that $|A\triangle B|$ is minimum. 
If $|A \triangle B|=2$ then by the second part of Lemma~\ref{lem:k-1}, 
there are $k(t-k+1)$ internally disjoint $AB$-paths in $F_k(G)$, implying $|\mathcal{C}|\geq k(t-k+1)$. 
Now assume that $|A \triangle B|=2r\geq 4$.
As in the proof of Theorem~\ref{thm:LowConnectivity}, since $r\geq 2$ and $n\geq\frac{1}{2}kt$ and $k^2-3k+4\geq0$. 
$$|\mathcal{C}|\geq 
r(n-k-r)+2r^2=
r(n-k+r)\geq 
2(n-k+2)\geq
kt -2k+4\geq 
k(t-k+1).$$
Thus $F_k(G)$ is $k(t-k+1)$-connected.
\end{proof}

The lower bound on the connectivity of $F_k(G)$ in Theorem~\ref{thm:BigConnectivity} is best possible.
For example, if $G$ is $t$-regular and contains a $k$-clique $X$, then $X$ has degree $k(t-k+1)$ in $F_k(G)$, 
implying $F_k(G)$ has connectivity at most $k(t-k+1)$. 
As a concrete example,  $G=K_{t+1}$ is $t$-connected, $t$-regular, and contains a $k$-clique. 
Thus the Johnson graph $J(t+1,k)\simeq F_k(K_{t+1})$ has connectivity at most $k(t-k+1)$.
In fact, the connectivity of $J(t+1,k)$ equals $k(t-k+1)$ \citep{MR1744235,diamjohn}. 
We conjecture the following generalization:

\begin{conj}
\label{conj:Connectivity}
If $G$ is a $t$-connected graph and $t\geq k$, then $F_k(G)$ is $k(t-k+1)$-connected. 
\end{conj}

Note that Conjecture~\ref{conj:Connectivity} with $k=2$
can be proved using the same method as the proof of
Theorem~\ref{thm:LowConnectivity} (since 
$|\mathcal{C}|\geq r(n-k-r) +2r^2=2(n-2-2)+8=2n> 2t>k(t-k+1)$). 

\section{Cliques}\label{sect:clique}

A \emph{clique} in a graph $G$ is a set of pairwise adjacent vertices in $G$. 
The \emph{clique number} $\omega(G)$ of $G$ is the maximum cardinality of a clique in $G$.
In this section we characterize the cliques in $F_k(G)$, and derive an exact formula for the clique-number of $F_k(G)$. 
These results are well known in the case of Johnson graphs \citep{choi}.

\begin{lemma}\label{lem:triangles}
Let $A,B,C$ be three pairwise adjacent vertices in $F_k(G)$. 
Then either $B\cap C \subset A$ or $A \subset B \cup C$  (but not both).
\end{lemma}

\begin{proof}
Suppose on the contrary that $B\cap C \not\subset A$ and $A \not\subset B \cup C$; 
that is, there are vertices $x\in (B\cap C)\setminus A$ and $a\in A\setminus (B \cup C)$. 
Since $A$ and $B$ are adjacent in $F_k(G)$ and $a\in A\setminus B$ and $x\in B\setminus A$, we have $A\triangle B=\{a,x\}$.
Similarly, $A\triangle C=\{a,x\}$.
Thus $B\cup C\cup\{a\}\setminus\{x\}\subseteq A$. 
Since $B$ and $C$ are adjacent in $F_k(G)$, we have $|B\cup C|=k+1$.
Thus $|A|\geq k+1$, which is the desired contradiction.
Thus $A\subset B \cup C$ or $B\cap C\subset A$. 

Now suppose that $A\subset B \cup C$ and $B\cap C\subset A$. 
Since $B$ and $C$ are adjacent, $|B\cap C|=k-1$ and $|B\cup C|=k+1$.
Since $|A|=k$, we have $A=B$ or $A=C$, which is the desired contradiction.
Thus $A \not\subset B \cup C$ or $B\cap C \not\subset A$.
\end{proof}

We now use Lemma~\ref{lem:triangles} to characterize cliques in $F_k(G)$.

\begin{theorem}\label{thm:charclique}
Let $X$ be a set of vertices in $F_k(G)$. 
Then $X$ is a clique of $F_k(G)$  if and only if there is a clique $K$ of $G$ and a set
$S\subseteq V(G)$, such that $K\cap S=\varnothing$ and either
\begin{enumerate}[(a)]
\item $X=\{ S \cup \{v\}  : v\in K \}$ and $|S|=k-1$, or
\item $X=\{(S \cup K)\setminus \{v\}  : v \in K\}$ and $|S|+|K|=k+1$.
\end{enumerate}
\end{theorem}

\begin{proof}
The ``if'' direction is immediate. To prove the ``only if'' direction, let $X$ be an arbitrary clique of $F_k(G)$. 

First suppose that $|X|=2$. Then $X=\{A,B\}$ for some edge $AB$ of $F_k(G)$. 
Let $S:=A\cap B$ and $K:=A\triangle B$. Then $X$ satisfies (a). In fact, it also satisfies (b). 






Now assume that $|X|=p\geq3$. Say $X=\{A_1,\dots,A_p\}$. 
For distinct $i,j\in[3,p]$, 
the two options given by Lemma~\ref{lem:triangles} 
for $A_1,A_2,A_i$ and $A_1,A_2,A_j$ are incompatible. 
That is, if say $A_i \subset A_1 \cup A_2$ but $A_1 \cap A_2\subset A_j$,
then $|A_i\triangle A_j|=4$, implying $A_i$ and $A_j$ are not adjacent in $F_k(G)$. 
Thus one of the following cases apply:

\begin{itemize}

\item $A_1\cap A_2\subset A_i$ for all $i\in[p]$: Let $S:=A_1\cap A_2$. 
Since $|S|=k-1$, each $A_i$ contains a vertex $v_i$ such that $A_i=S\cup\{v_i\}$. 
Thus $A_i\triangle A_j=\{v_i,v_j\}$  for distinct $i,j\in[p]$. 
Since $A_i$ and $A_j$ are adjacent in $F_k(G)$, $v_i$ and $v_j$ are adjacent in $G$. 
Thus $K:=\{v_i:i\in [p]\}$ is a clique in $G$, and $X=\{ S\cup\{v\}  : v\in K \}$.
Hence $X$ satisfies (a). 

\item $A_i\subset A_1\cup A_2$ for all $i\in[p]$: 
Since $|A\cup B|=k+1$, each $A_i$ contains a vertex $v_i$ such that $A_i=(A_1\cup A_2)\setminus\{v_i\}$. 
Thus $A_i\triangle A_j=\{v_i,v_j\}$  for distinct $i,j\in[p]$. 
Since $A_i$ and $A_j$ are adjacent in $F_k(G)$, $v_i$ and $v_j$ are adjacent in $G$. 
Thus $K:=\{v_i:i\in [p]\}$ is a clique in $G$.
Moreover, $X=\{ (S\cup K)\setminus\{v\}  : v\in K \}$ where $S:=(A_1\cup A_2)\setminus K$.
Hence $X$ satisfies (b). 
\end{itemize}
This completes the proof. Note that in both cases $S=\bigcap_i A_i$ and $K=\bigcup_iA_i\setminus S$. 
\end{proof}

We obtain the following formula for the clique-number of a token graph.

\begin{theorem}\label{thm:cliquenum}
$\omega(F_k(G)) = \min \{\omega(G), \max\{n-k+1,k+1\}\}$. 
\end{theorem}

\begin{proof}
We first prove the upper bound on $\omega(F_k(G))$.
Let $X$ be a clique in $F_k(G)$ with $\omega(F_k(G))$ vertices. 
Thus $X$ satisfies (a) or (b) in Theorem~\ref{thm:charclique}.
In case (a), $n\geq|S|+|K|=k-1+|X|$. 
In case (b), $|X|=|K|\leq |S|+|K|=k+1$.  
Thus $|X|\leq n-k+1$ or $|X|\leq k+1$, implying $|X|\leq\max\{n-k+1,k+1\}$. 
In both cases, $|X|=|K|\leq\omega(G)$. 
Therefore $\omega(F_k(G))=|X|\leq \min \{\omega(G), \max\{n-k+1,k+1\}\}$. 

We now prove the lower bound on $\omega(F_k(G))$.
Let $K$ be a clique in $G$ with $\omega(G)$ vertices. 
Consider the following two constructions of cliques in $F_k(G)$:
\begin{itemize}

\item Let $K'$ be a subset of $K$ with $\min\{\omega(G),n-k+1\}$ vertices.
Thus $|V(G)\setminus K'|\geq k-1$. Let $S$ be a subset of $V(G)\setminus K'$ with $k-1$ vertices.
Thus $\{S\cup\{v\}:v\in K'\}$ is a clique in $F_k(G)$ with $|K'|$ vertices.
Hence $\omega(F_k(G)) \geq \min\{\omega(G), n-k+1\}$. 


\item Let $K'$ be a subset of $K$ with $\min\{\omega(G),k+1\}$ vertices.
Since $n\geq k+1$, there is a subset $S$ of $V(G)\setminus K'$ with $(k+1)-|K'|$ vertices.
Thus $\{(S\cup K')\setminus\{v\}:v\in K'\}$ is a clique in $F_k(G)$ with $|K'|$ vertices.
Hence $\omega(F_k(G)) \geq \min\{\omega(G), k+1\}$. 

\end{itemize}

Therefore $\omega(F_k(G)) \geq \max\{\min\{\omega(G), n-k+1\}, \min\{\omega(G), k+1\}\}$, which equals $\min \{\omega(G), \max\{n-k+1,k+1\}\}$. 
\end{proof}


\begin{cor}
\label{cor:CliqueNumber}
Assuming $k\leq\frac{n}{2}$, we have $\omega(F_k(G)) = \min\{\omega(G),n-k+1\}$.
\end{cor}


For Johnson graphs, Corollary~\ref{cor:CliqueNumber} amounts to a
special case of the Erd\H{o}s-Ko-Rado theorem, which states that if
$0< t< k$ and  $\mathcal{F}$ is a family of $k$-subsets of an $n$-set
and $n\geq n_0(k,t)$ and the intersection of any two sets in
$\mathcal{F}$ has cardinality at least $t$, then
$|\mathcal{F}|\leq\binom{n-t}{k-t}$. \citet{wilson} proved this result
with $n_0(k,t)=(t+1)(k-t+1)$, which is best possible. Observe that a
clique in $J(n,k)$ is such a family $\mathcal{F}$ for $t=k-1$. In this
case, Wilson's Theorem states that $\omega(J(n,k))\leq n-k+1$ whenever $n\geq 2k$.

\section{Chromatic Number}\label{sect:chrom}


In this section we study the chromatic number of $F_k(G)$ in terms of the chromatic number of $G$. 
Our first result is an upper bound on $\chi(F_k(G))$.

\begin{theorem}
\label{thm:upperchrom}
$\chi(F_k(G))\le \chi(G)$.
\end{theorem}
\begin{proof}
Let $c:V(G)\rightarrow \{0,1, \dots ,\chi(G)-1\}$ be a coloring of $G$. 
To each vertex $A$ of $F_k(G)$, assign the color  $$c'(a):=\Big(\sum_{x \in A} c(x) \Big)\mod \chi(G)\enspace.$$ 
Let  $A$ and $B$ be two adjacent  vertices in $F_k(G)$.
Thus $A \triangle B =\{a,b\}$ for some edge $ab$ of  $G$. 
Suppose on the contrary that $c'(A)=c'(B)$.
Thus  $$\sum_{x \in A} c(x) \equiv \sum_{y \in B} c(y) \pmod{\chi(G)}\enspace.$$ 
Since $A \triangle B =\{a,b\}$, we have $c(a)\equiv c(b) \pmod{\chi(G)}$. 
Hence $c(a)=c(b)$, and $c$ is not a coloring of $G$. 
This contradiction proves that  $c'$  is a coloring of $F_k(G)$.
\end{proof}

Note that Theorem~\ref{thm:upperchrom} holds with equality whenever $\omega(G)=\chi(G)$ and $n\geq\omega(G)+k-1$, in which case $\chi(F_k(G))\geq\omega(F_k(G))\geq\omega(G)=\chi(G)$ by Theorem~\ref{thm:cliquenum}.

We now consider lower bounds on the chromatic number of token graphs.
By Theorem~\ref{thm:cliquenum}, we have
$\chi(F_k(G))\geq\omega(F_k(G))=\min \{\omega(G), \max\{n-k+1,k+1\}\}$. 
But we can obtain qualitatively stronger lower bounds in terms of $\chi(G)$ as follows.
First consider the case when $F_k(G)$ is bipartite.

\begin{prop}
If $F_k(G)$ is bipartite for some $k\geq1$, then $F_\ell(G)$ is bipartite for all $\ell\geq1$. 
\end{prop}

\begin{proof}
By Theorem~\ref{thm:upperchrom}, it suffices to prove that if $F_k(G)$ is bipartite then $G$ is bipartite. 
Equivalently, we prove that if $G$ is not bipartite then $F_k(G)$ is not bipartite. 
Suppose that $G$ is not bipartite.
Thus $G$ contains an odd cycle $C=(v_1,\dots,v_p)$.
First suppose that $p\geq k+1$.  Hence
\begin{align*}
& \{v_1,v_2,\dots,v_{k-2},v_{k-1},v_k\}
 \{v_1,v_2,\dots,v_{k-2},v_{k-1},v_{k+1}\}
 \{v_1,v_2,\dots,v_{k-2},v_{k-1},v_{k+2}\}\cdots\\
& \{v_1,v_2,\dots,v_{k-2},v_{k-1},v_{p}\}
 \{v_1,v_2,\dots,v_{k-2},v_{k},v_{p}\}
 \{v_1,v_2,\dots,v_{k-1},v_{k},v_{p}\}\cdots\\
& \{v_1,v_3,\dots,v_{k-1},v_{k},v_{p}\}
 \{v_2,v_3,\dots,v_{k-1},v_{k},v_{p}\}
\end{align*}
is a $p$-cycle in $F_k(G)$. Thus $F_k(G)$ is not bipartite.
Now assume that $p\leq k$. 
Let $A$ be a set of $k-p+1$ vertices in $V(G)\setminus C$ (which exist since $n\geq k+1$). 
Then $F_k(G,A)\simeq F_{p-1}(G-A)$ by \eqref{eqn:fixed}. 
Since $C$ is contained in $G-A$, by the above construction, 
there is an odd cycle in $F_{p-1}(G-A)$.
Thus there is an odd cycle in $F_k(G,A)$, which is a subgraph of $F_k(G)$. 
Thus $F_k(G)$ is not bipartite.
\end{proof}

We have the following general lower bound on $\chi(F_k(G))$. 

\begin{theorem} 
\label{thm:lowerchrom}
  $\chi(F_k(G)) \ge \frac{n-k+2}{n}\chi(G)-1.$
\end{theorem}
\begin{proof}
The result holds for $k=1$ since $F_1(G)\simeq G$.
Now assume that $k \ge 2$.
Let  $V_1, \dots V_{\chi(G)}$ be the colors classes in a coloring of $G$ with $\chi(G)$ colors. 
Assume that $|V_1| \ge \dots \ge |V_{\chi(G)}|$. 
Thus for each index $m$, 
\begin{equation}
\label{eqn:lowerchrom}
\sum_{i=1}^m |V_i| \geq \frac{mn}{\chi(G)}\enspace.
\end{equation}
Let $m$ be the minimum index such that  $\sum_{i=1}^m |V_i| \ge k-1$. Thus \eqref{eqn:lowerchrom} implies that
\begin{align}
\label{eqn:lowerchrom2}
\frac{(m-1)n}{\chi(G)} \leq \sum_{i=1}^{m-1} |V_i| \le k-2\enspace.
\end{align}
Let $X$ be a subset of $\bigcup_{i=1}^{m}V_i$ of cardinality $k-1$. 
Since $G[X]$ is $m$-colorable, 
\begin{align*}
\chi(G)  \leq \chi(G[X]) + \chi(G-X) \leq m + \chi(G-X) \enspace.
\end{align*}
By \eqref{eqn:fixed},  $G-X \simeq  F_1(G-X)\simeq F_k(G,X)$, which is a subgraph of $F_k(G)$. Thus
\begin{align*}
\chi(G) 
\leq m+\chi(F_k(G,X))
 \leq m+\chi(F_k(G)) \enspace.
\end{align*}
By \eqref{eqn:lowerchrom2},
\begin{align*}
\chi(G)  \leq\frac{k-2}{n}\chi(G)+1+\chi(F_k(G))\enspace.
\end{align*}
The result follows.
\end{proof}

Theorem~\ref{thm:lowerchrom} and \eqref{eqn:k,n-k} imply the following lower bound on $\chi(F_k(G))$ independent of $k$.

\begin{theorem}
\label{thm:CombinedLowerChrom}
$\chi(F_k(G)) \ge (\frac{1}{2}+\frac{2}{n})\chi(G)-1$ for all $k\geq1$.
\end{theorem}


Theorem~\ref{thm:CombinedLowerChrom} gives a lower bound of roughly $\frac12 \chi(G)$ on $\chi(F_k(G))$. 
However, the best upper bound example we know of is $\chi(F_k(G))\leq \chi(G)-2$, 
which is achieved for $G=K_n$ and $k=3$, for all $n >7$ and 
$n\equiv 1\pmod 6$ or $n\equiv 3\pmod 6$; see \citep{Lu,teirlinck,Ji-JCTA05}. 
In this case, an independent set in $J(n,3)$ is a Steiner triple system. 
\citet{etzion}  give some other values of $n$ and $k$ for which $\chi(J(n,k)) < n$.
These results suggest the following question, which is open even for Johnson graphs. 

\begin{open}
Does there exist a constant $c >0$ such that $\chi(F_k(G))\ge \chi(G)-c$ for every graph $G$ and
integer $k\geq1$?
\end{open}

\section{Hamiltonian Paths}\label{sect:ham}

In this section we study conditions for the existence or non-existence
of Hamiltonian paths in token graphs. 
First note that all Johnson graphs are Hamiltonian~\citep{Hans}.
Now consider the case when $G=P_n$, the path on $n$ vertices; see Figure~\ref{fig:TokenPath}. 
A Hamiltonian path in $F_k(P_n)$ would correspond to a Gray code of adjacent transpositions for the set of binary strings of length $n$ with $k$ ones. 
This Gray code exists if and only if $n$ is even and $k$ is odd; see \cite[p.~133]{ruskey} or \citep{Savage97}.  
Thus $F_k(P_n)$ contains a Hamiltonian path if and only if $n$ is even and $k$ is odd.
Hence: 

\begin{theorem}\label{thm:ham}
If a graph $G$ contains a Hamiltonian path and $n$ is even and $k$ is odd, then $F_k(G)$ contains a Hamiltonian path.
\end{theorem}

Note that the existence of a Hamiltonian cycle or path in $G$ does not imply that $F_k(G)$ contains a Hamiltonian cycle or path.
For example, $C_4$ is Hamiltonian, but $F_2(C_4)\simeq K_{2,4}$ does not even contain a Hamiltonian path.
More generally, if $G$ is bipartite and $\binom{n}{k}$ is odd\footnote{While most binomial coefficients are even, there are infinitely many non-trivial binomial coefficients that are odd; see \citep{MR1397551,MR0429714}.}, then $F_k(G)$ is bipartite by Theorem~\ref{thm:upperchrom}, but $F_k(G)$ is not Hamiltonian, since every bipartite Hamiltonian graph has even order. 
Even if $F_k(G)$ has  even order, it may not contain a Hamiltonian path. 
For example, let $V_1$ and $V_2$ be the color classes of $K_{m,m}$.  
Then $F_k(K_{m,m})$ is also bipartite by Theorem~\ref{thm:upperchrom}, and the color classes are 
\begin{align*}
W_1&=\{A \in V(F_k(K_{m,m} ) ):|A \cap V_1| \text{ is even} \} \text{ and }\\
W_2&=\{A \in V(F_k(K_{m,m} )): |A \cap V_1| \text{ is odd} \}\enspace.
\end{align*}
Thus, by an identity of  \citet{gould} (see \cite[p.~61]{Sprugnoli}), 
$$|W_1|-|W_2|=
	\sum_{i=0}^k (-1)^i \binom{m}{i}\binom{m}{k-i}
        =\begin{cases}0& \text{ if $k$ is odd}\\
(-1)^{k/2} \binom{m}{k/2} & \text{ if $k$ is even}
\enspace. \end{cases}$$
Hence for even $k$, $\big||W_1|-|W_2|\big|>2$ and therefore $F_k(K_{m,m})$ does not contain a Hamiltonian path. 
On the other hand, $F_k(K_{m,m})$  contains a Hamiltonian path for odd $k$ by Theorem~\ref{thm:ham}.

\section{Cartesian Product}\label{sect:prod}

The Cartesian product $G \square H$ of two graphs $G$ and $H$ is the graph with vertex set $V(G) \times V(H)$, where 
two vertices $(g,h)$ and $(g',h')$ are adjacent in $G \square H$ whenever $g=g'$ and $hh'\in E(H)$, or $h=h'$ and $gg'\in E(G)$. 
The Cartesian product of $m \ge 3$ graphs $G_1, \dots, G_m$ is defined recursively as $G_1 \square (G_2 \square \cdots \square G_m)$.
We now show that certain induced subgraphs in a token graph are in fact Cartesian products.

Let $H$ and $H'$ be two disjoint induced subgraphs of a graph $G$. 
Let $r$ and $s$ be integers such that $1\leq r \leq |V(H)|$ and $1\leq s \le |V(H')|$ and $r+s =k$.
Observe that the subgraph of $F_k(G)$ induced by all $k$-sets $A$ of $G$ such that $|A \cap V(H)|=r$  and $|A\cap V(H')|=s$ 
is isomorphic to $F_r(H) \square F_s(H')$. Thus $F_r(H) \square F_s(H')$ is an induced subgraph of $F_k(G)$.
We conclude:

\begin{theorem}\label{thm:cartmany}
If $H_1, \dots, H_m$ are pairwise disjoint induced subgraphs of a graph $G$,
then for all integers $s_1, \dots, s_m$ such that $1\leq s_i \le |V(H_i)|$ and $\sum s_i =k$, 
the graph $F_{s_1}(H_1) \square \cdots \square F_{s_m}(H_m)$ is an induced subgraph of $F_k(G)$.
\end{theorem}

In the case $k=2$, Theorem~\ref{thm:cartmany} has the following interpretation:

\begin{cor} \label{cor:cart2}
Let $H$ and $H'$ be two disjoint induced subgraphs of $G$. 
Then $H \square H'$ is an induced subgraph of $F_2(G)$.
\end{cor}

Corollary~\ref{cor:cart2} implies, for example, that the
$\lfloor\frac{n}{2}\rfloor\times\lceil\frac{n}{2}\rceil$ grid graph is
an induced subgraph of $F_2(P_n)$; see Figure~\ref{fig:TokenPath}. This shows that $F_2(G)$ can have unbounded treewidth even for trees $G$. Moreover, $F_2(G)$ can have unbounded clique minors even for trees $G$, since $F_2(K_{1,n})$ is isomorphic to $K_n$ with each edge subdivided once.

\section{Open Problems}\label{sect:end}

We now consider some open problems regarding $F_k(G)$ that are related to graph reconstruction. 
Does a given token graph uniquely determine the original graph?
We conjecture that this is indeed so.

\begin{conj}\label{conj:recon}
Let $G$ and $H$ be two graphs, such that $F_k(G) \simeq F_k(H)$ for some $k$. Then $G \simeq H$.
\end{conj}

This conjecture is related to the well known Reconstruction Conjecture; 
see \citep{BH-JGT77} for a survey. 
The \emph{deck} of a graph $G$ is the multiset of unlabeled
graphs $\{G-v:v\in V(G)\}$. The Reconstruction Conjecture states
that a graph is uniquely determined up to isomorphism by its deck.
Similarly, Conjecture~\ref{conj:recon} states that a graph is uniquely
determined up to isomorphism by one of its token graphs. 
Given that each element of the deck of $G$ is an induced subgraph of $F_2(G)$, 
it is possible that progress in this direction will shed light on the
Reconstruction Conjecture.


\medskip 
We conclude the paper with two definitions: 
For $r\in[k]$, let $F_{k,r}(G)$ be the graph with vertex set $\binom{V(G)}{k}$, 
where two vertices $A$ and $B$ in $F_{k,r}(G)$ are adjacent whenever 
$|A\triangle B|=2r$ and there is a perfect matching between  $A\setminus B$ and $B \setminus A$ in $G$. 
This graph is a generalization of the token graph since $F_k(G)\simeq F_{k,1}(G)$.
It is also a generalization of the Kneser graph $KG_{n,k}$, whose vertices are the $k$-subsets of an $n$-set, 
where two vertices $A$ and $B$ are adjacent whenever $A\cap B=\varnothing$. 
Observe that $KG_{n,k} \simeq F_{k,k}(K_n)$. 
Finally, let $F'_{k,r}(G)$ be the variant where instead we require that every edge is present between  $A\setminus B$ and $B \setminus A$.
Then again  $F_k(G)\simeq F'_{k,1}(G)$ and $KG_{n,k} \simeq F'_{k,k}(K_n)$. 
The study of $F_{k,r}(G)$ and $F'_{k,r}(G)$ is an open line of research.


\end{document}